\documentclass[12pt]{amsart}
\usepackage{amscd,amssymb}
\usepackage[graph,frame,poly,arc]{xy}  
\usepackage[plainpages,backref,urlcolor=blue]{hyperref}

\theoremstyle{plain}
\newtheorem{thm}[subsection]{Theorem}
\newtheorem{lem}[subsection]{Lemma}
\newtheorem{prop}[subsection]{Proposition}
\newtheorem{cor}[subsection]{Corollary}

\theoremstyle{definition}
\newtheorem{rk}[subsection]{Remark}

\newtheorem{conj}[subsection]{Conjecture}

\numberwithin{equation}{section}
\newcommand{\OO}{{\mathcal O}}

\newcommand{\A}{{\mathcal A}}

\newcommand{\al}{{\alpha}}
\newcommand{\be}{{\beta}}

\newcommand{\Z}{\mathbb{Z}}
\newcommand{\Q}{\mathbb{Q}}

\newcommand{\C}{\mathbb{C}}
\newcommand{\PP}{\mathbb{P}}

\begin{document}

\title [Line and rational curve arrangements, and Walther's inequality]
{Line and rational curve arrangements, and Walther's inequality}

\author[Alexandru Dimca]{Alexandru Dimca$^{1}$}
\address{Universit\'e C\^ ote d'Azur, CNRS, LJAD and INRIA, France }
\email{dimca@unice.fr}

\author[Gabriel Sticlaru]{Gabriel Sticlaru}
\address{Faculty of Mathematics and Informatics,
Ovidius University
Bd. Mamaia 124, 900527 Constanta,
Romania}
\email{gabrielsticlaru@yahoo.com }


\thanks{$^1$ This work has been partially supported by the French government, through the $\rm UCA^{\rm JEDI}$ Investments in the Future project managed by the National Research Agency (ANR) with the reference number ANR-15-IDEX-01}



\begin{abstract} There are two invariants associated to any line arrangement: the freeness defect $\nu(C)$
and an upper bound for it, denoted by $\nu'(C)$, coming from a recent result by Uli Walther. We show that $\nu'(C)$ is combinatorially determined, at least when the number of lines in $C$ is odd, while the same property is conjectural for $\nu(C)$.  In addition, we conjecture that the equality $\nu(C)=\nu'(C)$ holds if and only if the essential arrangement $C$ of $d$ lines has 
either a point of multiplicity $d-1$, or has only double and triple points. We prove both conjectures in some cases, in particular when the number of lines is at most 10. We also extend a result by H. Schenck on the Castenuovo-Mumford regularity of line arrangements to arrangements of possibly singular rational curves.

\end{abstract}
 
\maketitle


\section{Introduction}

Let $S=\C[x,y,z]$ be the graded polynomial ring in three variables $x,y,z$ with complex coefficients and let $C:f=0$ be a reduced curve of degree $d$ in the complex projective plane $\PP^2$. The minimal degree $mdr(f)$ of a Jacobian relation for the polynomial $f$ is the  smallest integer $m\geq 0$ such that there is a nontrivial relation
\begin{equation}
\label{rel_m}
 af_x+bf_y+cf_z=0
\end{equation}
among the partial derivatives $f_x, f_y$ and $f_z$ of $f$ with coefficients $a,b,c$ in $S_m$, the vector space of  homogeneous polynomials in $S$ of degree $m$. When $mdr(f)=0$, then $C$ is a union of lines passing through one point. We assume from now on that 
\begin{equation}
\label{eq1}
mdr(f)\geq 1.
\end{equation}
When $C$ is a line arrangement, this condition says exactly that $C$ is essential, i.e. not all the lines in $C$ pass through one point.
Denote by $\tau(C)$ the global Tjurina number of the curve $C$, which is the sum of the Tjurina numbers of the singular points of $C$. Recall that, if $(X,0)$ is an isolated plane curve singularity at the origin $0 \in \C^2$,
given by a local equation $g(u,v)=0$, where $g \in \OO_2$, the ring of analytic function germs at the origin, then the Tjurina number of $(X,0)$ is by defintion
\begin{equation}
\label{Tjur}
\tau(X,0)= \dim_{\C}\frac{\OO_2}{(g,g_u,g_v)},
\end{equation}
where $g_u$ and $g_v$ are the partial derivatives of $g$ with respect to the local coordinates $u$ and $v$.
We denote by $J_f$ the Jacobian ideal of $f$, i.e. the homogeneous ideal in $S$ spanned by $f_x,f_y,f_z$, and  by $M(f)=S/J_f$ the corresponding graded ring, called the Jacobian (or Milnor) algebra of $f$.
 Let $I_f$ denote the saturation of the ideal $J_f$ with respect to the maximal ideal ${\bf m}=(x,y,z)$ in $S$ and consider the local cohomology group 
 $$N(f)=I_f/J_f=H^0_{\bf m}(M(f)).$$
 The graded $S$-module  $N(f)$ satisfies a Lefschetz type property with respect to multiplication by generic linear forms, see  \cite{DPop}. This implies in particular the inequalities
\begin{equation}
\label{eqDP1} 
0 \leq n(f)_0 \leq n(f)_1 \leq ...\leq n(f)_{\lceil  T/2 \rceil} \geq n(f)_{\lceil T/2\rceil+1} \geq ...\geq n(f)_T \geq 0,
\end{equation}
where $T=3d-6$, $n(f)_k=\dim N(f)_k$ for any integer $k$, and for any real number $u$,  $\lceil     u    \rceil $ denotes the round up of $u$, namely the smallest integer $U$ such that $U \geq u$.
We set
$$\nu(C)=\max _j \{n(f)_j\},$$
and call $\nu(C)$ the {\it freeness defect}  of the curve $C$. It is known that a curve $C$ is free (resp. nearly free) if and only if $\nu(C)=0$ (resp. $\nu(C)=1$), see \cite{DStRIMS}.
When $d=2m$ is even, then the above implies that
$n(f)_{3m-3}=\nu(C).$
When $d=2m+1$ is odd, then the above and the self duality of the graded $S$-module $N(f)$, see \cite{Se, SW}, implies that
\begin{equation}
\label{eqDP2} 
 n(f)_{k} = n(f)_{T-k}, 
\end{equation}
for any integer $k$, and in particular
$n(f)_{3m-2}=n(f)_{3m-1}=\nu(C).$
The relation between the invariants $\nu(C)$ and $mdr(f)$ is given by the following result, see \cite[Theorem 1.2]{Drcc}.
\begin{thm}
\label{thmN}
Let $C:f=0$ be a reduced plane curve of degree $d$ and let $r=mdr(f)$.
Then the following hold.
\begin{enumerate}
\item If $r < d/2$, then
$$\nu(C)=(d-1)^2-r(d-1-r)-\tau(C).$$

\item If $r \geq (d-2)/2$, then 
$$\nu(C)= \lceil   \  \frac{3}{4}(d-1)^2 \    \rceil  -\tau (C).$$

\end{enumerate}
\end{thm}
Note that for $(d-2)/2 \leq r \leq (d-1)/2$, both formulas (1) and (2) above apply, and they give the same result for $\nu(C)$, as a direct simple check shows.

 Terao's conjecture  says that the freeness of a line arrangement is combinatorially determined. For more on Terao's conjecture and free hyperplane arrangements, we refer to \cite{DHA, Yo}. Note that Terao's conjecture would follow if the invariant $mdr(f)$ were combinatorially determined. 
But  examples due to G. Ziegler in \cite{Zi} show that 
the invariant $mdr(f)$ is not combinatorially determined, see also \cite[Example 4.3]{AD}. 

The following stronger version of H. Terao's conjecture was made in  \cite{Drcc}. 
\begin{conj}
\label{conj1}
Let $C:f=0$ be a line arrangement in $\PP^2$. Then the invariant $\nu(C)$ is combinatorially determined.
\end{conj}
Some cases where this conjecture holds are described in \cite[Proposition 3.5]{Drcc}.
Let $T\langle C\rangle$ denote the vector bundle of logarithmic vector fields along $C$, as defined and studied for instance in \cite{AD,Ang,DS14, FV}.
Let $(d_1^{L_0},d_2^{L_0})$, with  $d_1^{L_0} \leq d_2^{L_0}$, denote the splitting type of the vector bundle $T\langle C\rangle(-1)$ along a generic line $L_0$ in $\PP^2$. Then \cite[Theorem 1.1]{AD} says that 
$$(d-1)^2 -d_1^{L_0}d_2^{L_0}=\tau(C)+ \nu(C).$$
 Since the total Tjurina number $\tau(C)$ is clearly a combinatorial invariant for $C$ a line arrangement and since  $d_1^{L_0}+d_2^{L_0}=d-1$, Conjecture \ref{conj1} may be restated as follows.

\begin{conj}
\label{conj2}
Let $C:f=0$ be a line arrangement in $\PP^2$. Then the generic splitting type  $(d_1^{L_0},d_2^{L_0})$ of the vector bundle $T\langle C\rangle(-1)$ is combinatorially determined.
\end{conj}
This conjecture, in the form of a question, was made for the first time in  \cite[Question 7.12]{C+}.

The only general upper bound we know for the invariant $\nu(C)$ comes from  a deep result by U. Walther, see \cite[Theorem 4.3]{Wa}, bringing into the picture the monodromy of the Milnor fiber $F:f=1$ associated to the curve $C:f=0$. 
For an alternative proof of this result, see M. Saito  \cite[Theorem A.1, Appendix]{SAI}.
This upper bound is recalled in Theorem \ref{UWthm} below and is denoted by $\nu'(C)$.
When $C$ is a rational cuspidal curve, this bound $\nu'(C)$ is very sharp, and this fact was used in \cite{DStRIMS, DStRIMS2} to show that most, and very likely all, rational cuspidal curves are either free or nearly free.
In this note we investigate how sharp is this upper bound in the case of line arrangements. 
 The following conjecture, suggested by the analysis of many examples, looks quite surprising to us, since it says that the equality $\nu(C)=\nu'(C)$ holds only in two extreme cases.

\begin{conj}
\label{conj3}
Let $C:f=0$ be an essential arrangement of $d\geq 3$ lines in $\PP^2$, and let $m(C)$ be the maximal multiplicity of an intersection point on $C$.
Then $\nu(C)=\nu'(C)$ if and only if either $m(C)=d-1$, or $m(C) \leq 3$.
\end{conj}
In the second section we recall basic facts on the spectrum of a line arrangement and state U. Walther's result mentioned above. As a direct application, we extend a result by H. Schenck on the Castenuovo-Mumford regularity of line arrangements to arrangements of possibly singular rational curves, see Theorem \ref{HSthm}.

In the third section we show that Conjecture \ref{conj3} holds, if we add an extracondition on the line arrangement $C:f=0$, namely that the invariant $mdr(f)$
satisfies the inequality $mdr(f) \geq (d-2)/2$, see Theorem \ref{thmA}. This condition is verified by all line arrangement $C$ having only double and triple points, and by many other line arrangements.
In the fourth section we show that Conjecture \ref{conj3} holds inside the class
of free arrangements with $mdr(f) \geq 2$, see Corollary \ref{corFREE}. And in the final section we show that all the three Conjectures above hold for arrangements having at most 10 lines, see Theorem \ref{thmD}.

\medskip

The first author thanks AROMATH team at INRIA Sophia-Antipolis for excellent working conditions, and in particular Laurent Bus\'e for stimulating discussions.

We would like to thank the referee for the very careful reading of our manuscript and for his very useful suggestions to improve the presentation.

\section{The spectrum of a line arrangement and Walther's inequality} 

Let $\A$ be an arrangement of $d$ planes in $\C^{3}$,  given by a reduced equation
$$ f(x,y,z)=0.$$ We assume $\A$ to be central, i.e. any plane in $\A$ passes through the origin $0 \in \C^3$. This implies that $f$ is a homogeneous polynomial of degree $d$, and the curve $C:f=0$ in $\PP^2$ is the associated line arrangement to $\A$.
Consider the corresponding global Milnor fiber $F$, defined by $f(x,y,z)-1=0$ in $\C^{3}$, with monodromy action $h:F \to F$, 
$$h(x,y,z)=\exp(2\pi i/d)\cdot (x,y,z).$$
 In studying the cohomology $H^*(F,\Q)$ of the Milnor fiber or the monodromy action
$h^*: H^*(F,\Q) \to H^*(F,\Q)$, we can, without any loss of generality,  suppose that the arrangement $\A$ is essential, and we do this in the sequel.
For basic facts on mixed Hodge structures in this setting we refer to  \cite{DHA}.
For $\be$ a $d$-th root of unity, let $H^j(F,\C)_{\beta}$ denote the corresponding monodromy eigenspace.

The interplay between the monodromy $h^*: H^*(F,\Q) \to H^*(F,\Q)$ and the mixed Hodge structure (MHS) on $H^*(F,\Q)$ is reflected by  the spectrum of $\A$  defined as 
\begin{equation}
\label{eqSpec}
{\rm Sp}(\A)=\sum_{\al \in \Q}m_{\al}t^{\al},
\end{equation}
with
\begin{equation}
\label{eqSpec2}
m_{\al}=\sum_j(-1)^{j}\dim {\rm Gr}_F^p\tilde H^j(F,\C)_{\beta},
\end{equation}
where $p=\lfloor 3-\al \rfloor$,
$\beta=\exp(-2\pi i\al)$ and $\tilde H^j(F,\C)$ is the reduced cohomology of the Milnor fiber $F$. Here $\lfloor x \rfloor$ denotes the integral part of the rational number $x$. Recall that each cohomology group $\tilde H^j(F,\C)$ admits a decreasing Hodge filtration
\begin{equation}
\label{eqSpec3}
\tilde H^j(F,\C)=F^0\tilde H^j(F,\C) \supset F^1\tilde H^j(F,\C) \supset \ldots \supset F^j\tilde H^j(F,\C) \supset F^{j+1}\tilde H^j(F,\C)=0, 
\end{equation}
and we set as usual
\begin{equation}
\label{eqSpec4}
{\rm Gr}_F^p\tilde H^j(F,\C)  =\frac{F^p\tilde H^j(F,\C)}{F^{p+1}\tilde H^j(F,\C)}. 
\end{equation}
Since the monodromy $h$ is a morphism of algebraic varieties, there is an induced linear mapping $h^*:{\rm Gr}_F^p\tilde H^j(F,\C) \to {\rm Gr}_F^p\tilde H^j(F,\C)$, and one considers the corresponding eigenspaces
${\rm Gr}_F^p\tilde H^j(F,\C)_{\beta}$ in the formula \eqref{eqSpec2} above.
The rational number $\al$ is called a spectral number for the plane arrangement $\A$, if $m_{\al} \ne 0$ in ${\rm Sp}(\A)$. The corresponding integer $m_{\al}$ is called the multiplicity of $\al$.
Note that the cohomology groups $H^m(F,\Q)_1=\ker(h^*-Id)$ and $H^{m}(F,\Q)_{-1}=\ker (h^*+Id)$ are mixed Hodge substructures in $H^m(F,\Q)$, as $(h^*-Id)$ and $(h^*+Id)$ are MHS morphisms. 
The paper \cite{BS} gives the following very simple formulas for the coefficients
$m_{\alpha}$.
Let $\nu_j$ be the number of points in the projective line arrangement $C$ associated to $\A$, of multiplicity $j,\,\,j\geq 3.$ 

\begin{thm} \label{jumpcoef} Assume that the arrangement $C$ is essential. Then
the multiplicity $m_{\alpha} = 0$ if either $\alpha \notin (0,3)$ or $\alpha d \notin \Z,\,$ where $d=|\A|$ is the number of lines in $C$. For a rational number $\alpha = \frac{e}{d} \in\, ]0,1]\,$ with $e \in [1,d] \cap \Z$, one has the following.
\begin{center}
$m_{\alpha} = \displaystyle{\binom {e-1} {2}} - \displaystyle{\sum_j \nu_j \binom{ \lceil ej/d \rceil -1} {2}} ,$\\
$m_{\alpha +1} = (e-1)(d-e-1) - \displaystyle{\sum_j} \nu_j(\lceil ej /d \rceil -1)(j- \lceil ej /d \rceil ),$ \\
$m_{\alpha+2} = \displaystyle{\binom {d-e-1} {2} }- \displaystyle{\sum_j \nu_j \binom{ j-\lceil ej /d \rceil } {2}} - \delta_{e,d},$
\end{center}
where $\lceil \beta \rceil:= \min\{k\in \Z\,|\, k\geq \beta\},\,$ and $\delta_{e,d}=1$ if $e=d$ and $0$ otherwise.
\end{thm}

The key result of U. Walther in \cite[Theorem 4.3]{Wa}  yields the following inequality.
\begin{thm}
\label{UWthm} For any reduced plane curve $C:f=0$, one has
$$\dim N(f)_{2d-2-j} \leq \dim Gr_F^1H^2(F, \C)_{\lambda},$$
for $j=1,2,...,d$ and  $\lambda= \exp(2\pi i (d+1-j)/d)=\exp(-2\pi i (j-1)/d)$. 
In particular, this inequality gives {\rm in the middle range}, that is for  $j_0=\lceil \frac{d+1}{2} \rceil$, the following
$$\nu(C)= \dim N(f)_{2d-2-j_0} \leq \nu'(C):= \dim Gr_F^1H^2(F, \C)_{\lambda_0},$$
where $\lambda_0= \exp(2\pi i (d+1-j_0)/d)=\exp(-2\pi i (j_0-1)/d)$. 
\end{thm}

\subsection{The case $d$ even, $C$ line arrangement }
Suppose first that $d$ is even, say $d=2d'$ and take $j=d'+1$ in Theorem \ref{UWthm}. 
Then $\dim N(f)_{2d-2-j}=\nu(C)$ and 
${\lambda}=-1$, in other words we get the following.
\begin{equation}
\label{M1}
\nu'(C) = \dim Gr_F^1H^2(F, \C)_{-1}.
\end{equation}
Using equation \eqref{eqSpec2}, we get
$$m_{3/2}=-\dim {\rm Gr}_F^1 H^1(F,\C)_{-1}+\dim {\rm Gr}_F^1 H^2(F,\C)_{-1}= -\dim {\rm Gr}_F^1 H^1(F,\C)_{-1}+\nu'(C) .$$
Since $H^1(F,\C)_{-1}$ is a pure Hodge structure of weight 1, see \cite[Theorem 7.7]{DHA}, it follows that
$$\dim {\rm Gr}_F^1 H^1(F,\C)_{-1}=\frac{1}{2} \dim H^1(F,\C)_{-1}.$$
Use next the second equation in Theorem \ref{jumpcoef} for $\al=1/2$, that is for $e=d'$, and get
$$m_{3/2}=(d'-1)^2-\sum_j\nu_j ( \lceil j/2 \rceil -1)( j-\lceil j/2 \rceil ).$$
Any point $p\in C$ of multiplicity $k_p\geq 2$ gives rise to a $k_p$-ordinary multiple point, whose local defining equation in suitable local coordinates is a homogeneous polynomial of degree $k_p$. The Tjurina number $\tau(C,p)$ of such a singular point  is given by
\begin{equation}
\label{Milnor1}
\tau(C,p)=(k_p-1)^2.
\end{equation}
For $j=2j_1$ even, we have
$$( \lceil j/2 \rceil -1)( j-\lceil j/2 \rceil )=(j_1-1)j_1$$
and the Tjurina number of a point $p$ of multiplicity $j=2j_1$ is given by the formula \eqref{Milnor1}, namely 
$$\tau(C,p)=(j-1)^2=4j_1 ^2-4j_1+1.$$
Therefore, one has
\begin{equation}
\label{Milnor2}
( \lceil j/2 \rceil -1)( j-\lceil j/2 \rceil )=\frac{1}{4}(\tau(C,p)-1).
\end{equation}
Similarly, for $j=2j_1+1$ odd, we have
$$( \lceil j/2 \rceil -1)( j-\lceil j/2 \rceil )=j_1^2$$
and the Tjurina number of a point $p$ of multiplicity $j=2j_1+1$ is 
$$\tau(C,p)=(j-1)^2=4j_1^2.$$
Therefore, in this case one has
\begin{equation}
\label{Milnor3}
( \lceil j/2 \rceil -1)( j-\lceil j/2 \rceil )=\frac{1}{4}\tau(C,p).
\end{equation}
If we sum the equalities in \eqref{Milnor2} and \eqref{Milnor3} for all the multiple points of the arrangement $C$, we get the following
\begin{equation}
\label{M1.1}
m_{3/2}=\frac{1}{4} \left( (d-2)^2-\tau(C) + \sum_{j \text{ even }}\nu_j \right).
\end{equation}
Finally we get, for $d$ even, the equality
\begin{equation}
\label{M1.2}
\nu'(C)=  \frac{1}{2} \dim H^1(F,\C)_{-1} + \frac{1}{4} \left( (d-2)^2-\tau(C) + \sum_{j \text{ even }}\nu_j \right).
\end{equation}

\subsection{The case $d$ odd, $C$ line arrangement}
Suppose now that $d$ is odd, say $d=2d'+1$ and take $j=d'+1$. As above, we get
\begin{equation}
\label{M2}
\nu'(C) = \dim  {\rm Gr}_F^1H^2(F, \C)_{\eta}.
\end{equation}
 for  ${\eta= \exp(-2\pi i d'/(2d'+1))= \exp(2\pi i (d'+1)/(2d'+1))   }$.
Using equation \eqref{eqSpec2}, we get
$$m_{\alpha +1} =-\dim {\rm Gr}_F^1 H^1(F,\C)_{\eta}+\dim {\rm Gr}_F^1 H^2(F,\C)_{\eta},$$
where $\al=d'/(2d'+1).$ Since $\eta$ has order $d$, it follows that $H^1(F,\C)_{\eta}=0$
and $\dim  {\rm Gr}_F^1H^2(F, \C)_{\eta}=\dim H^{1,1}(F,\C)_{\eta}$, see \cite[Proposition 7.4]{DHA}.
Use Theorem \ref{jumpcoef} for $e=d'$ and get
$$m_{\alpha +1}=d'(d'-1)-\sum_j\nu_j ( \lceil jd'/(2d'+1) \rceil -1)( j-\lceil jd'/(2d'+1) \rceil ).$$
For a point $p$ of multiplicity $j=2j_1$ even, we have exactly as above
$$( \lceil jd'/(2d'+1) \rceil -1)( j-\lceil jd'/(2d'+1) \rceil )  =(j_1-1)j_1= \frac{1}{4}(\tau(C,p)-1).$$
Similarly, for a point $p$ of multiplicity $j=2j_1+1$ odd, we have
$$( \lceil j/2 \rceil -1)( j-\lceil j/2 \rceil )=j_1^2=\frac{1}{4}\tau(C,p),$$
since $j_1<d'$.
These formulas imply, by summation over all the multiple points $p$, the following.
\begin{equation}
\label{M2.1}
m_{\alpha +1} =\frac{1}{4} \left( (d-1)(d-3)-\tau(C) + \sum_{j \text{ even }}\nu_j  \right).
\end{equation}
Finally we get, for $d$ odd, the equality
\begin{equation}
\label{M2.2}
\nu'(C)= \frac{1}{4} \left( (d-1)(d-3)-\tau(C) + \sum_{j \text{ even }}\nu_j  \right).
\end{equation}

\begin{rk}
\label{rkComb}
Note that for $d$ odd, the invariant $\nu'(C)$ is clearly determined by the combinatorics, see formula \eqref{M2.2}. For $d$ even, the same is true, in view of  formula \eqref{M1.2}, if the multiplicity of the eigenvalue $-1$ for the monodromy acting on $H^1(F,\C)$ is determined by combinatorics, a fact that is often conjectured, see \cite{PS,S2}.
\end{rk}

\subsection{An application of Walther's inequality to rational curve arrangements} 

We conclude this section by giving an application of Walther's inequality.
Recall the definition of the {\it stability threshold} 
$$st(f)=\min \{q~~:~~\dim M(f)_k=\tau(C) \text{ for all } k \geq q\},$$
for any reduced plane curve $C:f=0$.
\begin{thm}
\label{HSthm}
Let $C:f=0$ be a plane curve of degree $d$, such that any irreducible component of $C$ is a rational curve. Then $n(f)_k=0$ for $k\leq d-3$ or $k \geq 2d-3$. Moreover,
$st(f) \leq 2d-4$, with equality when $C$ is in addition nodal.
\end{thm}
 H. Schenck has proven this result in  \cite[Corollary 3.5]{HS}, but only for line arrangements, and he has phrased it in terms of the Castelnuovo-Mumford regularity $reg(f)$ of the Milnor algebra $M(f)$. The relation between $st(f)$ and $reg(f)$ is clearly explained in \cite[Theorem 3.4]{DIM}: one has 
$reg(f)=st(f)$ if $C$ is a free curve, and $reg(f)=st(f)-1$ otherwise.
Note that the proof of H. Schenck is quite different from ours, as it uses induction on the number of lines.
\proof
Apply Theorem \ref{UWthm} for $j=1$, and note that in this case $\lambda=1$.
It is known that $H^2(F,\C)_1=H^2(U,\C)$, where $U$ is the complement $\PP^2 \setminus C$, see for instance \cite[Equation (5.2.20), p.147]{DSTH}. Since all the components of $C$ are rational curves it follows from \cite[Theorem 2.7]{Ab}
that 
$$ \dim Gr_F^1H^2(F, \C)_1= \dim Gr_F^1H^2(U, \C)=0.$$
Then Theorem \ref{UWthm} implies that $n(f)_{2d-3}=0$. The vanishings of the integers $n(f)_k$ claimed in Theorem \ref{HSthm} follows from the inequalities \eqref{eqDP1} and the equalities \eqref{eqDP2}. The claim about $st(f)$ follows using the formula
$$n(f)_k = \dim M(f)_k + \dim M(f)_{T-k}-\dim M(f_F)_k -\tau(C),$$
where $f_F=x^d+y^d+z^d$, see \cite[Equation (2.8)]{DStRIMS}. The claim of the equality in the case of nodal curves follows from
\cite[Theorem 1.3]{DStCamb}.
\endproof

\section{Walther's inequality for  line arrangements with double and triple points} 
In this section we compute the invariants $\nu(C)$ and $\nu'(C)$ for various classes of line arrangements $C$. To refer to certain line arrangements in $\PP^2$, we recall the following notation from \cite{DIM}. We say that a line arrangement $C$ of $d\geq 4$ lines is of type $L(d,m)$ if there is a single point of multiplicity $m \geq 3$ and all the  other intersection points of $C$  are double points. We say that $C$ is of type $\hat L(m_1,m_2)$, with $3 \leq m_1 \leq m_2$ if there is a line $L$ in $C$ containing two points, say $p_1$ of multiplicity $m_1\geq 3$ and $p_2$ of multiplicity $m_2 \geq 3$, and such that all the other lines in $C$ pass through either $p_1$ or $p_2$.
The monomial line arrangement $\A(m,m,3)$ for $m\geq 2$ is given by the equation
$$C: f=(x^m-y^m)(x^m-z^m)(y^m-z^m)=0.$$
We recall the following result. 

\begin{prop}
\label{propA}
Let $C:f=0$ be a  line arrangement of $d$ lines in $\PP^2$. Then one has the following.

\begin{enumerate}
\item $mdr(f)=1$ if and only if $d=3$ and $C$ is a triangle, or $d \geq 4$ and $C$ is of type $L(d,d-1)$. Any such arrangement is free.

\item Any arrangement $C$ of type $L(d,d-2)$ for $d\geq 5$ is nearly free, with $mdr(f)=2$.

\item Any arrangement $\hat L(m_1,m_2)$ is free with $mdr(f)=m_1-1$.

\item An arrangement $C:f=0$ has $mdr(f)=2$ if and only if $C$ is either of type $L(d,d-2)$, or of type 
$\hat L(3,m_2)$, or linearly equivalent to the monomial line arrangement $\A(2,2,3)$.

\item The  line arrangement $\A(m,m,3)$ is free with  $mdr(f)=m+1$.

\end{enumerate}

\end{prop}

For the claims (1) and (2) we refer to \cite[Proposition 4.7]{DIM}, for (3) see 
\cite[Example 4.10]{DIM}, for (4) we refer to \cite{To} or \cite[Theorem 4.11]{DIM}, and finally, for (5) see \cite[Example 8.6]{DHA}.
In this section we prove the following main result.

\begin{thm}
\label{thmA}
Let $C:f=0$ be a line arrangement of $d\geq 3 $ lines in $\PP^2$.

\begin{enumerate}
\item If the line arrangement $C:f=0$ has
 only double and triple points, then $mdr(f) \geq (d-2)/2$ and hence
 $$\nu(C)= \lceil   \  \frac{3}{4}(d-1)^2 \    \rceil  -\tau (C).$$
 In particular, Conjecture \ref{conj1} and Conjecture \ref{conj2} hold in this situation.
\item For a line arrangement $C:f=0$ satisfying  $mdr(f) \geq (d-2)/2$,
one has
$$\nu(C)=\nu'(C),$$
i.e. Walther's inequality in the middle range is an equality, if and only if the line arrangement $C$ has
 only double and triple points.
 
 \end{enumerate}
\end{thm}

\proof
First we address the claim (1).
For a generic line arrangement $C:f=0$, i.e. an arrangement having only points of multiplicity 2, we know that $r=mdr(f)=d-2$, see \cite[Theorem 4.1]{DStEdin}.
If $C$ has only double and triple points, then for $d=4$, we have in addition to the generic arrangement the possibility $L(4,3)$, when $r=1$.
For $d=5$, we have two additional cases besides the generic arrangement, namely $L(5,3)$, 
and $\hat L(3,3)$, both with $r=2$.
For a line arrangement having only double and triple points, we can use the inequality in \cite[Example 2.2 (iii)]{DS14} and conclude that
$$r=mdr(f) \geq 2d/3-2.$$
For $d \geq 6$, we get $2d/3-2 \geq (d-2)/2$, and hence in all cases we can use the formula (2) in Theorem \ref{thmN} and get
$$\nu(C)= \lceil   \  \frac{3}{4}(d-1)^2 \    \rceil  -\tau (C).$$
Written down explicitly, this means that for $d=2d'$ even, one has
$\nu(C)=3(d')^2-3d'+1-\tau(C)$, while for $d=2d'+1$ odd, one has
$\nu(C)=3(d')^2-\tau(C).$

Now we address the claim (2). Let $C:f=0$ be a line arrangement satisfying $mdr(f) \geq (d-2)/2$. The corresponding invariant $\nu(C)$ is then exactly as in the claim (1).  
On the other hand, the corresponding integers $\nu_j\geq 0 $ satisfy the following two relations
		$$\sum_{k\geq 2} \nu_k{k \choose 2}={d \choose 2}
		\text{ and } \sum_{k\geq 2} \nu_k(k-1)^2=\tau(C).$$
		The first equality is well known, see for instance \cite[Exercise 2.8]{DHA}, and the second follows from the fact that the Tjurina number of a point of multiplicity $k$ equals $(k-1)^2$, as stated in \eqref{Milnor1}.
By eliminating $\nu_3$, we get
$$\nu_2=2d(d-1)-3\tau(C)+ \sum_{k\geq 4}\nu_k(k^2-4k+3).$$
Note that $k^2-4k+3>0$ for $k \geq 4$, hence the last sum, which we denote by $\Sigma$, is strictly positive if $C$ has points of multiplicity $>3$.
For $d=2d'$ even and $C$ a line arrangement having only double and triple points, note that $H^1(F,\C)_{-1}=0$, see \cite[Corollary 5.4]{DHA}.
The formula \eqref{M1.2} yields $\nu'(C)=3(d')^2-3d'+1-\tau(C)$ for  a line arrangement having only double and triple points
and
$$\nu'(C)=  \frac{1}{2} \dim H^1(F,\C)_{-1} + 3(d')^2-3d'+1-\tau(C)+\frac{1}{4}\Sigma$$
in general. 
Similarly, for $d=2d'+1$ odd,
the formula \eqref{M2.2} yields $\nu'(C)=3(d')^2-\tau(C)$ for  a line arrangement having only double and triple points
and
$$\nu'(C)=   3(d')^2-\tau(C)+\frac{1}{4}\Sigma$$
in general. 
These formulas imply the claim (2).
\endproof


\section{Walther's inequality for  free line arrangements} 

 Let now $C:f=0$ be a free line arrangement with exponents $(d_1,d_2)$, which means that $T\langle C\rangle$, the vector bundle of logarithmic vector fields along $C$, splits as a direct sum, namely
 $T\langle C\rangle(-1)=\OO_{\PP^2}(-d_1) \oplus \OO_{\PP^2}(-d_2)$. It is known that $d_1+d_2=d-1$ and 
 also $d_1=mdr(f)$, if we assume, as we'll do in the sequel, that $d_1 \leq d_2$.
For a free curve $C$ with exponents $(d_1, d_2)$, one has  $\nu(C)=0$, see \cite[Proposition 8.2]{DHA}, and hence  Theorem \ref{thmN} implies that  
$$\tau(C)=(d_1+d_2)^2-d_1d_2.$$
For free line arrangements of even degree $d$ we have the following.
\begin{thm}
\label{thmB}
Let $C:f=0$ be a free line arrangement in $\PP^2$ with exponents $(d_1,d_2)$, where $d=d_1+d_2+1$ is even. Then one has
$$0=\nu(C) <\nu'(C),$$
i.e. Walther's inequality  in the middle range is strict, unless we are in one 
of the following  cases, when
$\nu'(C)=0$.

\begin{enumerate}
\item $d_1=1$ and any $d$;

\item $d_1=2$, $d=6$ and  $C$ is linearly equivalent to the monomial arrangement $\A(2,2,3)$.

\end{enumerate}

\end{thm}

\proof

For $d=2d'$, the formula \eqref{M1.2} yields
$$\nu'(C)=  \frac{1}{2} \dim H^1(F,\C)_{-1} + \frac{1}{4} \left( d_1d_2-2(d_1+d_2)+1 + \sum_{j \text{ even }}\nu_j \right).$$
We list now the cases when $\nu'(C)=0$. To do this, note that $\nu'(C)=0$ implies 
\begin{equation}
\label{Triv1}
d_1d_2-2(d_1+d_2)+1=(d_1-2)(d_2-2)-3 \leq 0.
\end{equation}

\begin{itemize}
		\item[i)] If $d_1=1$, then $C$ consists of $(d-1)$ lines passing through a common point, and an additional secant, see \cite{DIM}. Then $d_2=d-2$, $\nu_2=d-1$ and $H^1(F,\C)_{-1}=0$, see \cite[Corollary 5.4]{DHA}. It follows that $\nu'(C)=0$.
		\item[ii)] If $d_1=2$, then $d \geq 6$ and the possibilities for $C$ are recalled in Proposition \ref{propA}: the corresponding intersection lattice should be $L(d,d-2)$, or $\hat L(3,d-2)$ or $C$ is the monomial arrangement $\A(2,2,3)$. For the intersection lattice $L(d,d-2)$, the arrangements are nearly free, but not free, see \cite[Proposition 4.7]{DIM}. For the intersection lattice $\hat L(3,d-2)$, we get $\nu_2=2(d-3)$ and $\nu_{d-2}=1$. Here $d \geq 6$, hence $\nu'(C)>0$ in this case. For the remaining case $\A(2,2,3)$ we get $\nu'(C)=0$, since $d_1=2$, $d_3=3$, $\nu_2=3$ and $H^1(F,\C)_{-1}=0$, for instance  by \cite[Corollary 5.4]{DHA}.
		\item[iii)] If $d_1=3$, the only case which may yield $\nu'(C)=0$ is when $d_2=4$, in view of \eqref{Triv1} and the fact that $d_4$ must be even. But any arrangement of $d=8$ lines, which is free with $d_1=3$ satisfies $\sum_{j \text{ even }}\nu_j>1$. Indeed, using \cite[Theorem 1.2]{Dcurvearr}, it follows that the maximal multiplicity of a point in such an arrangement is at most $5$. The integers $\nu_j\geq 0 $ satisfy the following two relations
		$$\nu_2+3\nu_3+6\nu_4+10\nu_5={8 \choose 2}=28$$
		$$\nu_2+4\nu_3+9\nu_4+16\nu_5=\tau(C)=37.$$
		We show first that $\nu_2=0$ leads to a contradiction.
		If we multiply the first equation by 4, the second by 3, make the difference and we assume $\nu_2=0$, we get $-3\nu_4-8\nu_5=1$, which is a contradiction. Assume now that $\nu_4=0$. Then we get $\nu_3+6\nu_5=9$, which gives rise to two cases: $A: \  \ (\nu_2,\nu_3,\nu_5)=(9,3,1)$ and $B: \  \  (\nu_2,\nu_3,\nu_5)=(1,9,0)$. To show that the case $B$ is impossible, note that these data imply that there is a line $L$ containing only triple points of the arrangement $C$. But then the intersection number of this line $L$ with the curve formed by the remaining 7 lines should be an even number, a contradiction since this intersection number is 7.
This shows that  $\nu'(C)>0$ when $d_1=3$.

\item[iv)] When $d_1 \geq 4$, we have $d_2 \geq d_1 \geq 4$ and hence
the condition \eqref{Triv1} is not fulfilled. Hence $\nu'(C)>0$ when $d_1 \geq 4$.
	\end{itemize}
	
\endproof
For free line arrangements of odd degree $d$ we have the following.
\begin{thm}
\label{thmC}
Let $C:f=0$ be a free line arrangement in $\PP^2$, with exponents $(d_1,d_2)$, where $d=d_1+d_2+1$ is odd. Then one has
$$0=\nu(C) <\nu'(C),$$
i.e. Walther's inequality  in the middle range is strict, unless we are in one 
of the following  cases, when
$\nu'(C)=0$.
\begin{enumerate}
\item $d_1=1$ and any $d$;

\item $d_1=2$, $d=5$ and the intersection lattice of $C$ is of type $\hat L(3,3)$;

\item $d_1=3$, $d=7$ and the arrangement is linearly equivalent to an arrangement obtained from the monomial line arrangement $\A'=\A(2,2,3)$, by adding a line determined by two double points in $\A'$.

\item $d_1=4$, $d=9$ and the arrangement is linearly equivalent to the monomial line arrangement $\A(3,3,3)$

\end{enumerate}

\end{thm}

\proof

For $d=2d'+1$, the formula \eqref{M2.2} yields
$$\nu'(C)=  \frac{1}{4} \left( d_1d_2-2(d_1+d_2) + \sum_{j \text{ even }}\nu_j \right).$$
We list now the cases when $\nu'(C)=0$. First note that $\nu'(C)=0$ implies 
\begin{equation}
\label{Triv2}
d_1d_2-2(d_1+d_2)=(d_1-2)(d_2-2)-4 \leq 0.
\end{equation}

\begin{itemize}
		\item[i)] If $d_1=1$, then $C$ consists of $(d-1)$ lines passing through a common point, and an additional secant, see \cite{DIM}. Then $d_2=d-2$, $\nu_2=d-1$ and $\nu_{d-1}=1$. It follows that $\nu'(C)=0$.
		\item[ii)] If $d_1=2$, we proceed as above. For the intersection lattice $\hat L(3,d-2)$, we get $\nu_2=2(d-3)$. Here $d \geq 5$, hence $\nu'(C)=0$ if and only if $d=5$. 
		\item[iii)]If $d_1=3$, the only cases which may yield $\nu'(C)=0$ are $d_2=3$ and $d_2=5$. Indeed, $d_2$ must be an odd number such that the pair $(d_1=3, d_2)$ satisfies the condition \eqref{Triv2}.
		Consider first the case $d_2=3$, i.e. $d=7$.
		Using \cite[Theorem 1.2]{Dcurvearr}, it follows that the maximal multiplicity of a point in such an arrangement is at most $4$. The integers $\nu_j\geq 0 $ satisfy the following two relations
		$$\nu_2+3\nu_3+6\nu_4={7 \choose 2}=21$$
		$$\nu_2+4\nu_3+9\nu_4=\tau(C)=27.$$
		If we multiply the first equation by 4, the second by 3, make the difference, we get $\nu_2-3\nu_4=3$. If  $\nu_4=0$, then $\nu_2=3$, $\nu_3=6$ and we get the arrangement described in point (3). If $\nu_4>0$, then we get $\nu'(C)>0$ as well.

We consider now the case $d_2=5$, i.e. $d=9$ and show that in this situation $\nu'(C)>0$. 
To do this it is enough to show that  $\sum_{j \text{ even }}\nu_j>1$. As above, we see that the maximal multiplicity of a point in such an arrangement is at most $6$ and we get the following relations. 
$$\nu_2+3\nu_3+6\nu_4+10\nu_5+ 15\nu_6={9 \choose 2}=36$$
		$$\nu_2+4\nu_3+9\nu_4+16\nu_5+25\nu_6=\tau(C)=49.$$
Exactly as above we get the relation
$$-\nu_2+3\nu_4+8\nu_5+15\nu_6=3.$$		
If $\nu_5 >0$ or $\nu_6>0$, the claim is clear since we get $\nu_2\geq 5$. If we assume
$\nu_5 =\nu_6=0$, then only the case $\nu_4=1$, $\nu_2=0$ need further discussion.
In this situation we get $\nu_3=10$, and these data contradict the 
 Hirzebruch inequality (which holds provided that $\nu_d=\nu_{d-1}=0$), see \cite[Remarks added in proof]{H},
\begin{equation}
\label{eq:hirzebruch}
\nu_2 +\frac{3}{4}\nu_3 \geq d+ \Sigma_{k \geq 5}(k-4)\nu_k.
\end{equation}

\item[iv)] If $d_1=4$, then the only possibility to get $\nu'(C)=0$ is $d_2=4$, in view of \eqref{Triv2}.
Hence $d=9$, and
in addition $\nu_2=\nu_4=\nu_6= \nu_8=0$. The relations
$$3\nu_3+10\nu_5={9 \choose 2}=36 \text{ and }
		4\nu_3+16\nu_5=\tau(C)=48,$$
imply $\nu_5=0$, and hence $C$ is linearly equivalent to 	$\A(3,3,3)$.

\item[v)] When $d_1 \geq 5$, we have $d_2 \geq d_1 \geq 5$ and hence
the condition \eqref{Triv2} is not fulfilled. Hence $\nu'(C)>0$ when $d_1 \geq 5$.			
			\end{itemize}
	
\endproof

\begin{cor}
\label{corFREE}
Conjecture \ref{conj3} holds for any free, essential line arrangement in $\PP^2$.
\end{cor}

\section{Conjecture \ref{conj3} holds for arrangements of $d \leq 10$ lines} 

We start with the following.

\begin{lem}
\label{lemL}
Let $C:f=0$ be a  line arrangement of type $L(d,d-2)$ for $d \geq 6$.
Then $\nu'(C)>\nu(C)=1$.
\end{lem}
\proof
By Proposition \ref{propA} (2) we know that $C$ is nearly free, and hence $\nu(C)=1$.
To compute $\nu'(C)$, we consider first the case $d=2d'$ even. Then one has $\nu_2=2(d-2)+1$, $\nu_{d-2}=1$ and $\tau(C)=d^2-4d+6$. The formula \ref{M1.2} implies
$$\nu'(C) \geq \frac{d-2}{2} \geq 2>1= \nu(C).$$
Similarly, for $d=2d'+1$ odd we get using formula \ref{M2.2}
$$\nu'(C) = \frac{d-3}{2} \geq 2>1= \nu(C).$$ 
\endproof

Using this Lemma, we can prove the following.

\begin{thm}
\label{thmD}
Let $C:f=0$ be an arrangement of $d$ lines in $\PP^2$. If  $3 \leq d\leq 10$, then Conjecture \ref{conj1}, Conjecture \ref{conj2} and 
Conjecture \ref{conj3} hold for $C$.
\end{thm}
\proof
First recall that the freeness of a line arrangement $C$ is determined by the combinatorics for $d \leq 10$, see \cite{ACKN, DIM13,FV} and the references given there. For such arrangements, the generic splitting type of the rank two vector bundle $T\langle C\rangle(-1)$ coincides with the exponents, and they are determined by the combinatorics. Moreover we have $\nu(C)=0$.
In view of this remark and of Corollary \ref{corFREE}, we may assume in this proof that $C$ is not free. 

Let $m(C)$ be the maximal multiplicity of an intersection  point in $C$.
Recall that  \cite[Theorem 1.2]{Dcurvearr} says that, if the arrangement $C$ is not free, then either

\medskip

\noindent ($A$) $mdr(f)=d-m(C)$, or

\medskip

\noindent ($A'$) $m(C) \leq mdr(f) \leq d-m(C)-1$. In particular, in this case  $m(C) \leq (d-1)/2$.

\medskip

We split the proof into three parts. In each case, to prove that the  Conjecture \ref{conj1} and  Conjecture \ref{conj2} hold, we have to show that, if there are two line arrangements $C:f=0$ and $C':f'=0$ with the same combinatorics, with at least one of them in case $A'$ above,
then $\nu(C)=\nu(C')$. Note that examples due to G. Ziegler in \cite{Zi} show one may have $mdr(f) \ne mdr(f')$ in such a situation, see also \cite[Example 4.3]{AD}. 

\medskip

\noindent {\bf The cases $d \leq 8$.}

Note that in the case $A'$ we have $2m(C) \leq 7$, and hence $m(C) \leq 3$.
Hence, if we have two arrangements $C$ and $C'$ as above, then
$m(C)=m(C')\leq 3$ and the claim follows from Theorem \ref{thmA}.

Hence from now on suppose that $C$ is in case $A$, and prove only Conjecture \ref{conj3}.
If $mdr(f) \leq 2$, then the claim follows from Proposition \ref{propA} and Lemma \ref{lemL}. On the other hand, if $mdr(f) \geq (d-2)/2$, 
then the claim follows from Theorem \ref{thmN} (2) and Theorem \ref{thmA} (2). Since $(d-2)/2 \leq 3$, these two cases cover all the possibilities when $3 \leq d \leq 8$.

\medskip

\noindent {\bf The case $d=9$.}
Note that in the case $A'$ we have $2m(C) \leq 8$, and hence $m(C) \leq 4$. If $m(C) \leq 3$, the claim follows
Theorem \ref{thmA}. If $m(C)=4$, then in both cases $A$ and $A'$ we get $mdr(f)\geq (d-2)/2=3.5$, and the claim follows from Theorem \ref{thmN} (2) and Theorem \ref{thmA} (2).

Hence from now on suppose that $C$ is in case $A$, and prove only Conjecture \ref{conj3}.
The cases $mdr(f) \leq 2$ and $mdr(f) \geq (d-2)/2=3.5$ can be treated as in the previous case.
Therefore we have to consider only the case $mdr(f)=3$ and $m(C)=6$. 

First note that, if there are at least two points $O_1$ and $O_2$ of multiplicity $6$ in a line arrangement $C$, then the number $d$ of lines in $C$ must be at least $6+6-1=11$. Indeed, at most one line in the two sets of 6 lines passing through $O_1$ and respectively $O_2$ can be in common.
So let $O$ be the unique point of multiplicity 6 in our arrangement $C$ of 9 lines, and denote by $C'$ the union of the 6 lines in $C$ passing through $O$. Then $C$ is obtained from $C'$ by taking the union with an arrangement of 3 lines, denoted by $C_3$. There are two cases to discuss.

\begin{itemize}
		\item[i)] $C_3$ has a point of multiplicity 3, say $P$, and this point is situated on a line in $C'$. Then the arrangement $C$ is of type $\hat L(4,6)$, and hence is free. This case is not possible by our assumption. 
		
		\item[ii)] In all the other cases $C$ has only points of multiplicity $\leq 3$ except the point $O$, and $\nu_3 \in \{0,1,2,3\}$.
		The relation 
		$$\nu_2+3\nu_3+ 15\nu_6={9 \choose 2}$$
		implies $\nu_2=21-3\nu_3$. This yields
		$$\tau(C)=(21-3\nu_3)+4\nu_3+25=46+\nu_3$$
		and then 
		$$\nu(C)=64-3(9-3-1)-(46+\nu_3)=3-\nu_3.$$
		On the other hand, a direct computation using formula \ref{M2.2} yields 
		$$\nu'(C)=6-\nu_3 >\nu(C).$$

\end{itemize}

\noindent {\bf The case $d=10$.}

For the same reasons as in the case $d=9$, we have to consider only the case $mdr(f)=3$ and $m(C)=7$, and prove only Conjecture \ref{conj3}. 

 First note that, if there are at least two points $O_1$ and $O_2$ of multiplicity $m(C)=7$ in a line arrangement $C$, then the number $d$ of lines in $C$ must be at least $7+7-1=13$. Indeed, at most one line in the two sets of 7 lines passing through $O_1$ and respectively $O_2$ can be in common.
So let $O$ be the unique point of multiplicity 7 in $C$ and denote by $C'$ the union of the 7 lines in $C$ passing through $O$. Then $C$ is obtained from $C'$ by taking the union with an arrangement of 3 lines, denoted by $C_3$. 

There are two cases to discuss.

\begin{itemize}
		\item[i)] $C_3$ has a point of multiplicity 3, say $P$, and this point is situated on a line in $C'$. Then the arrangement $C$ is of type $\hat L(4,7)$, and hence is free. This case is not possible by our assumption. 
		
		\item[ii)] In all the other cases $C$ has only points of multiplicity $\leq 3$ except the point $O$, and $\nu_3 \in \{0,1,2,3\}$.
		The relation 
		$$\nu_2+3\nu_3+ 21\nu_7={10 \choose 2}$$
		implies $\nu_2=24-3\nu_3$. This yields
		$$\tau(C)=(24-3\nu_3)+4\nu_3+36=60+\nu_3$$
		and then 
		$$\nu(C)=81-3(10-3-1)-(60+\nu_3)=3-\nu_3.$$
		On the other hand, a direct computation using formula \ref{M1.2} yields 
		$$\nu'(C) \geq 7-\nu_3 >\nu(C).$$

\end{itemize}

\endproof

\begin{rk}
\label{rk11lines}
To  check our Conjectures for line arrangements of $d=11$ lines, one can try to use the same approach as above. The difficult case to discuss is when $m(C)=4$, since both values $mdr(f)=7$ and $mdr(f)=4$ seem to be possible, and they lead to distinct values for $\nu(C)$, namely $76-\tau(C)$ and respectively $75-\tau(C)$, by Theorem \ref{thmN}.
Hence not even the proof of Conjecture \ref{conj1} and  Conjecture \ref{conj2} can be obtained in this way, and a detailed analysis of this situation is necessary.

\end{rk}

\end{document}